\documentclass[a4paper,11pt]{article}
\usepackage{mathrsfs}
\usepackage{pgfplots}
\usepackage{amsfonts}
\usepackage{graphicx}
\usepackage{tikz}
\usepackage{float}
\usepackage{color}
\usepackage{ifpdf}
\usepackage{epstopdf}
\usepackage{transparent}
\usepackage{amsmath,amsfonts,amssymb,amscd}
\usepackage{indentfirst,graphics,epsfig}
\usepackage[comma,numbers,square,sort&compress]{natbib}

\input{epsf}
\makeatletter \@addtoreset{figure}{section} \makeatother
\makeatletter
\long\def\@makecaption#1#2{%
   \vskip 10\p@
   \setbox\@tempboxa\hbox{{#1}\ \ #2}%
   \ifdim \wd\@tempboxa >\hsize
       {#1}\ \ #2\par
   \else
       \hbox to\hsize{\hfil\box\@tempboxa\hfil}%
   \fi}
\makeatother
\newtheorem{thm}{Theorem}
\newtheorem{prop}{Proposition}

\newtheorem{lem}{Lemma}

\makeatletter \@addtoreset{equation}{section}

\newcommand{\qed}{{\hfill\rule{3pt}{7pt}}}
\def\pf{\noindent {\it Proof.} }

\begin{document}
\rule{0cm}{3.5cm}
\begin{center} {\Large \bf  On the general sombor index of connected unicyclic graphs with given diameter }
\end{center}

\pagestyle{empty} \vskip 2mm
\begin{center}
{
  {\small Xipeng Hu, Lingping Zhong}\\

  {\small Department of Mathematics}\\
  {\small Nanjing University of Aeronautics and Astronautics}\\

  {\small Key Laboratory of Mathematical Modelling }\\
  {\small and High Performance Computing of Air Vehicles}
  {\small (NUAA), MIIT, Nanjing, 211106, China}\\
  {\small Email addresses: huxipeng@nuaa.edu.cn, zhong@nuaa.edu.cn }\\
   }
\end{center}

\begin{abstract}
 \noindent
 The general sombor index of $G$ is defined as $SO_{\alpha}(G)= \sum_{uv\in G}\left(d^2_{G}+d^2_{G}\right)^{\alpha}$. We obtain upper bounds on the general sombor index $SO_{\alpha}(G)$ for unicyclic graph $G$ with given diameter and $0<\alpha<1$. The extremal graph was also given.
 \\[2mm]
 \end{abstract}

\section{Introduction}

Let $G$ be a graph with vertex set $V(G)$ and edge set $E(G)$.
$N_{G}(v)$ is the set consist of vertices which are adjacent to $v\in V(G)$. The degree of $v\in V(G)$ is $|N_{G}(v)|$.
A pendant vertex is a vertex of degree one. $VP(G)$ is the set consist of all pendant vertices in $G$. $N_{G}(v)$ of $v\in V(G)$ is the set of vertices adjacent to $v$ in $G$. 
 We denote the diameter path of $G$ by $diam(G)$.
 Denote a unicyclic graph with $n$ vertices and $diam(G)=d$ by $\mathcal{U}_{n,d}$. The girth $g(G)$ of a unicyclic graph $G$ is the length of the cycle in $G$. For other standard graph-theoretical notions, these can be found in \cite{J.A}.

Recently, a new vertex degree based topological index called Sombor index was introduced by Gutman, defined as\cite{IG21}
$$ SO(G)=\sum_{uv\in E(G)}\left(d^2_{G}(u)+d^2_{G}(v)\right)^{\frac{1}{2}}.$$
Mathematical results about sombor index was given in \cite{HE22}\cite{LZ22}\cite{SA21}\cite{TZ21}\cite{TZar}\cite{XJ22}.
In this paper, we define the general sombor index as follows:
$$ SO_{\alpha}(G)=\sum_{uv\in E(G)}\left(d^2_{G}(u)+d^2_{G}(v)\right)^{\alpha}.$$
If $\alpha=\frac{1}{2}$, it's sombor index.  If $\alpha=1$, we get the forgotten index $F(G)$.
$$ F(G) = \sum_{uv\in E(G)}\left(d^2_{G}(u)+d^2_{G}(v)\right)$$

In \cite{QB19}\cite{MK21}\cite{MK22}\cite{QB19}\cite{YM21}\cite{ZH18}\cite{HE22}\cite{LZ22}, determine the extremal value and its corresponding graph of connected graphs(unicyclic graphs, trees) with given diameter. Motivated by \cite{MK21,MK22}, we  determine the maximum general sombor index of unicyclic graph with given diameter.

\section{Preliminary results}

\begin{lem}\label{lem1}
Let $1\leq x $, $1\leq y$. For $\alpha>0$, we have
$$\left(x^2+y^2\right)^\alpha \leq \left(1+(x-1+y)^2\right)^\alpha.$$
\end{lem}

\pf \begin{align*}
\left(x^2+y^2\right)^\alpha&=\left((1+(x-1))^2+y^2\right)^\alpha =\left(1+2(x-1)+(x-1)^2+y^2\right)^\alpha\\
&< \left(1+(x-1)^2+2(x-1)y+y^2\right)^\alpha=\left(1+\left((x-1)+y\right)^2\right)^\alpha.
\end{align*}

\begin{figure}[H]
 \centering
  \includegraphics[width=0.5\textwidth]{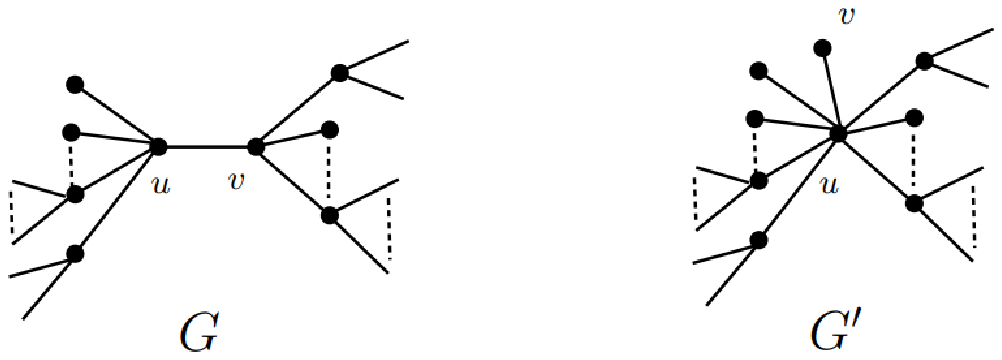}

   \caption{$G \stackrel{uv}{\longrightarrow} G'$ transformation}
    \label{fig:1}
     \end{figure}

\begin{lem}\label{lem2}
Let $G$ be a graph,  $uv\in E(G)$,  and $d_{G}(u)\geq 2$, $d_{G}(v)\geq 2$, 
 Let $G'$ be a graph obtained from $G$ by  relocating the component connected with $v$ except $u$ to $u$, let $v$ be a pendant vertex, we denoted this transformation by $G \stackrel{uv}{\longrightarrow}G'$ transformation(see Fig~\ref{fig:1}), then $SO_{\alpha}(G)<SO_{\alpha}(G')$
\end{lem}

\begin{figure}[H]
 \centering
  \includegraphics[width=0.5\textwidth]{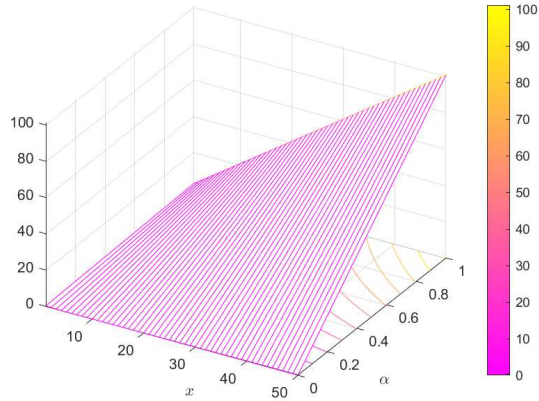}

   \caption{the function $g(x)$ for $0<\alpha<1$ and $x\geq 1$}
    \label{fig:4}
     \end{figure}

\pf  By Lemma~\ref{lem1},
\begin{align*}
SO_{\alpha}(G')-SO_{\alpha}(G)=&\sum_{xu\in E(G)\setminus \left\{uv\right\}}\left(\left(d_{G}^2(x)+\left(d_{G}(u)+d_{G}(v)-1\right)^2\right)^\alpha-\left(d_{G}^2(x)+d_{G}^2(u)\right)^\alpha\right)\\
&+
\sum_{yv\in E(G)\setminus \left\{uv\right\}}\left(\left( d_{G}^2(y)+(d_{G}(u)+d_{G}(v)-1)^2 \right)^\alpha-\left(d_{G}^2(y)+d_{G}^2(v)\right)^\alpha  \right)\\
&+
\left(1+(d_{G}(u)+d_{G}(v)-1)^2\right)^\alpha-\left( d_{G}^2(u) + d_{G}^2(v) \right)^\alpha>0,
\end{align*}
so we get $SO_{\alpha}(G)<SO_{\alpha}(G')$.
\begin{lem}\label{lem5}
$f(x)=(x^2+9)^\alpha- (x^2+4)^\alpha$ is a decreasing function for $x>0$ and $0<\alpha<1$.
\end{lem}
\pf$$ f'(x) = 2\alpha x\left(\left( x^2+9\right)^{\alpha-1}- \left( x^2+4\right)^{\alpha-1}\right)<0,$$ for $x>0$ and $0<\alpha<1$. So the conclusion holds.
\begin{lem}\label{lem6}
$f_{1}(x) = (x-1)\left(\left(x+1\right)^2+1\right)^{\alpha}+2\left(\left(x+1\right)^2+4\right)^{\alpha}-(x-1)\left(x^2+1\right)^{\alpha}- \left(x^2+4\right)^{\alpha} $  and $f_{2}(x) = (x-1)\left(\left(x+1\right)^2+1\right)^{\alpha}+2\left(\left(x+1\right)^2+4\right)^{\alpha}-(x-1)\left(x^2+1\right)^{\alpha}- \left(x^2+9\right)^{\alpha} $   are strictly increasing  for $x\geq 1$ and $0<\alpha<1$.
\end{lem}
\pf\begin{align*}
f_{1}'(x)=& 2\alpha(x-1)(x+1)\left( x^2+2x+2\right)^{\alpha-1}+4\alpha(x+1)\left(x^2+2x+5\right)^{\alpha-1}-2\alpha x(x-1)\left(x^2+1\right)^{\alpha-1}\\
&-2\alpha x\left(x^2+4\right)^{\alpha-1}+\left( x^2+2x+2\right)^{\alpha} - \left( x^2+1\right)^{\alpha}\\
\geq& 2\alpha\Big[(x-1)(x+1)\left( x^2+2x+2\right)^{\alpha-1}+2(x+1)\left(x^2+2x+5\right)^{\alpha-1}\\
&- x(x-1)\left(x^2+1\right)^{\alpha-1}- x\left(x^2+4\right)^{\alpha-1}  \Big].
\end{align*}
Let $g(x)=(x-1)(x+1)\left( x^2+2x+2\right)^{\alpha-1}+2(x+1)\left(x^2+2x+5\right)^{\alpha-1}- x(x-1)\left(x^2+1\right)^{\alpha-1}- x\left(x^2+4\right)^{\alpha-1}$. $g(x)> 0$ by computer drawing(see Fig~\ref{fig:4}). Hence, $f_{1}'(x)>0$  and $f_{1}(x)$ is a strictly increasing function for $0<\alpha<1$. Similarly, we can prove $f_{2}(x)$ is strictly increasing.
\begin{lem}\label{lem7}
$f(x)= (x-2)\left(x^2+1\right)^{\alpha}+ 2\left(x^2+4\right)^{\alpha}-2\left(\left(x-1\right)^2+4\right)^{\alpha}-(x-3)\left(\left(x-1\right)^2+1\right)^{\alpha}$ is a strictly increasing function for $x\geq 3$ and $0<\alpha<1$.
\end{lem}
\pf
\begin{align*}
f'(x)=& 2\alpha(x-2)x\left(x^2+1\right)^{\alpha} +4\alpha x\left(x^2+4 \right)^{\alpha-1}+\left(x^2+1\right)^{\alpha}
-2\alpha(x-3)(x-1)\left( \left( x-1\right)^2+1\right)^{\alpha-1}\\
 &- \left(\left( x-1\right)^2+1\right)^{\alpha}-4\alpha(x-1)\left( \left(x-1\right)^2+4\right)^{\alpha}\\
 \geq &2\alpha \Big[(x-2)x\left(x^2+1\right)^{\alpha}+2x\left(x^2+4 \right)^{\alpha-1}-(x-3)(x-1)\left( \left( x-1\right)^2+1\right)^{\alpha-1}\\
 &-2(x-1)\left( \left(x-1\right)^2+4\right)^{\alpha}\Big].
\end{align*}
Let $h(x)=(x-2)x\left(x^2+1\right)^{\alpha}+2x\left(x^2+4 \right)^{\alpha-1}-(x-3)(x-1)\left( \left( x-1\right)^2+1\right)^{\alpha-1}-2(x-1)\left( \left(x-1\right)^2+4\right)^{\alpha}$. $h(x)>0$ by computer drawing. So $f'(x)>0$ and $f(x)$ is a strictly increasing function.
\begin{figure}[H]
 \centering
  \includegraphics[width=0.5\textwidth]{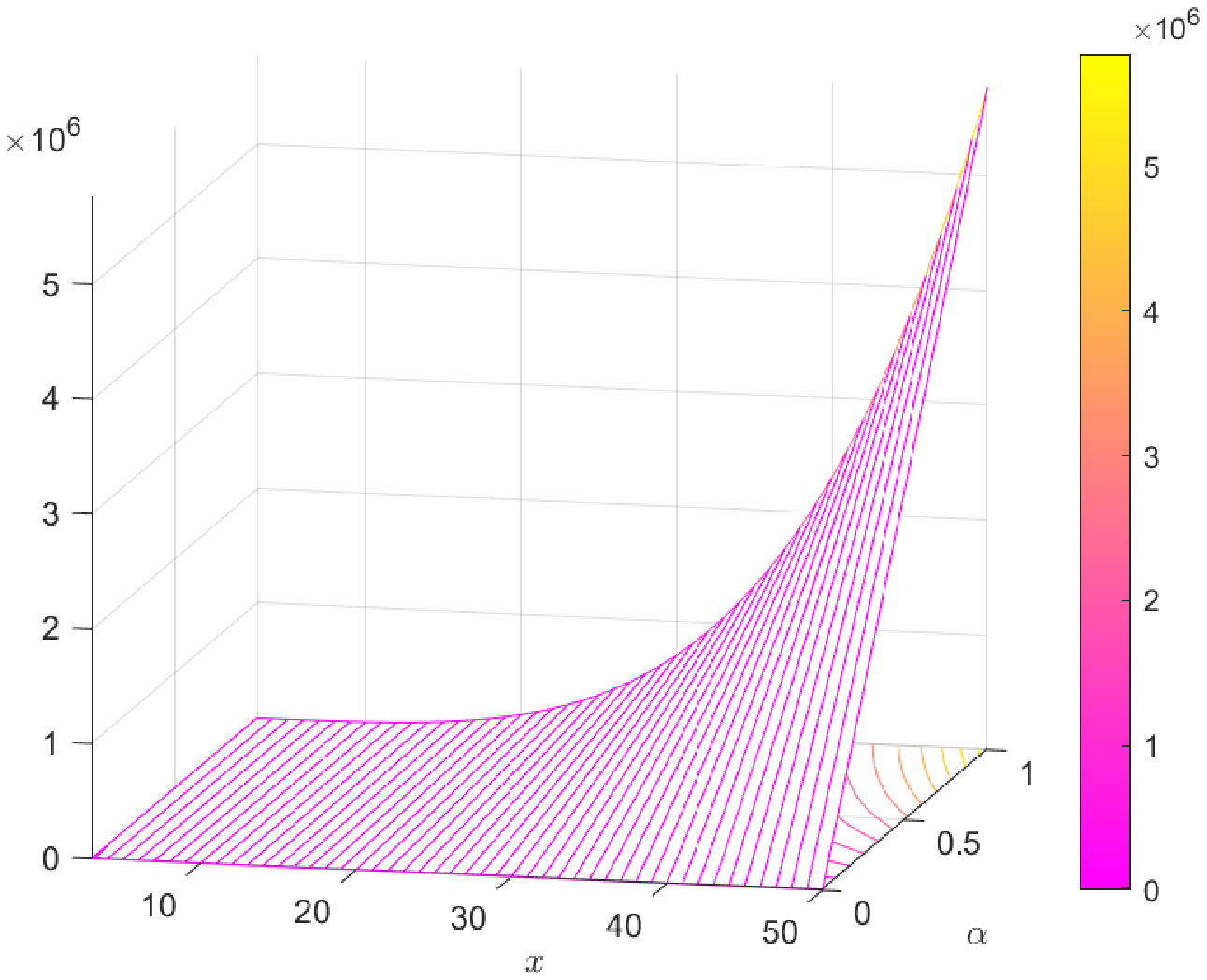}

   \caption{the function $h(x)$ for $0<\alpha<1$ and $x\geq 3$}
    \label{fig:5}
     \end{figure}

\section{Main results}
\begin{lem}\label{lem3}
Let $n\geq6$, the cycle $C_{n} \notin \mathcal{U}^{max}_{n,\lfloor \frac{n}{2}\rfloor}$.
\end{lem}
\pf Let $C_{n} = u_{0}u_{1}u_{2}\cdots u_{n-1}u_{0}$.

If $n$ is even, we use transformation $C_{n} \stackrel{u_{n-1}u_{0}}{\longrightarrow} G^{'}$, then $SO_{\alpha}(C_{n})<SO_{\alpha}(G')$, by Lemma~\ref{lem2}.

If $n$ is odd, let $G' = C_{n}+u_{1}u_{n-2}-u_{0}u_{n-1}$, then
$$SO_{\alpha}(C_{n})-SO_{\alpha}(G')=2(8^\alpha-13^\alpha)+8^\alpha-18^\alpha+2(8^\alpha-10^\alpha)<0.
$$
So  $C_{n} \notin \mathcal{U}^{max}_{n,\lfloor \frac{n}{2}\rfloor}$.\qed
\begin{lem}\label{lem4}
Let $G^*\in \mathcal{U}^{max}_{n,d}$, $3\leq d \leq n-2$ and $0<\alpha<1$, $C$ and
$P=v_{0}v_{1}\cdots v_{d}$ are the cycle and diametral path in $G^{*}$,  respectively. If $VP(G^*)\subseteq N_{G^*}(v_{1})\cup N_{G^*}(v_{d-1}) $, then $|V(C)\cap V(P)|>1$.
\end{lem}
\pf If $|V(C)\cap V(P)|\leq 1$,  we consider the following cases.

\textbf{Case 1.} If $V(C)\cap V(P) = \emptyset$, then there is a path, say $v_{i}u_{1}\cdots u_{l}$ and $l\geq 1$, connecting $P$ and $C$ in $G^*$.  Using the transformation $G^{*} \stackrel{u_{l}u_{l-1}}{\longrightarrow} G^{*'}$,  we get the graph $G^{*'}\in \mathcal{U}_{n,d}$ , by Lemma~\ref{lem2},
$SO_{\alpha}(G^{*})<SO_{\alpha}(G^{*'})$,  which is a contradiction.

\textbf{Case 2.} If $V(C)\cap V(P) = 1$. Suppose  $V(C)\cap V(P) = v_{i}$,   $g(G^{*})=k$.  Denote the cycle of $G^{*}$ by $C=v_{i}w_{1}w_{2}\cdots w_{k-1}v_{i}$.

\textbf{Subcase 2.1.} If $k\geq 4$,  we use the transformation $G^{*}\stackrel{v_{i}w_{1}}{\longrightarrow} G^{*'}$,  get the graph $G^{*'}$  which $SO_{\alpha}(G^{*})<SO_{\alpha}(G^{*'})$  by Lemma~\ref{lem2}, a contradiction.

\textbf{Subcase 2.2.} If $k = 3$, we know $d_{G^*}(v_{i+1})\geq 2$, $d_{G^*}(v_{i})\geq 4$ and $d_{G^*}(v_{i+2})\geq 1$. Denote the cycle of $G^*$ by $C= v_{i}w_{1}w_{2}v_{i}$.     Let $G^{*'}=G^*+v_{i+1}w_{1}-w_{1}w{2}$, then we have
\begin{align*}
SO_{\alpha}(G^*)-SO_{\alpha}(G^{*'})=& 8^{\alpha} - \left( 4+\left(d_{G^*}(v_{i+1})+1\right)^2\right)^{\alpha}
+\left(d_{G^*}^2(v_{i+1})+4\right)^{\alpha}-\left(d_{G^*}^2(v_{i+1})+1\right)^{\alpha}\\
&+\left(d_{G^*}^2(v_{i})+d_{G^*}^2(v_{i+1})\right)^{\alpha}-\left(d_{G^*}^2(v_{i})+\left(d_{G^*}^2(v_{i+1})+1\right)^2\right)^{\alpha}\\
&+\left(d_{G^*}^2(v_{i+1})+d_{G^*}^2(v_{i+2})\right)^{\alpha}-\left(\left(d_{G^*}^2(v_{i+1})+1\right)^2+d_{G^*}^2(v_{i+2})\right)^{\alpha}\\
\leq&8^{\alpha} - \left( 4+\left(d_{G^*}(v_{i+1})+1\right)^2\right)^{\alpha}+\left(d_{G*}^2(v_{i+1})+4\right)^{\alpha}-\left(d_{G*}^2(v_{i+1})+1\right)^{\alpha}\\
\leq&8^{\alpha} - 13^{\alpha} +20^{\alpha}- 17^{\alpha},
\end{align*}
we know  $8^{\alpha} - 13^{\alpha} +20^{\alpha}- 17^{\alpha}<0$ when $0<\alpha<1.90056$. Hence, $SO_{\alpha}(G^*)<SO_{\alpha}(G^{*'})$, which is a contradiction. So we complete the proof.\qed

\begin{figure}[H]
 \centering
  \includegraphics[width=0.5\textwidth]{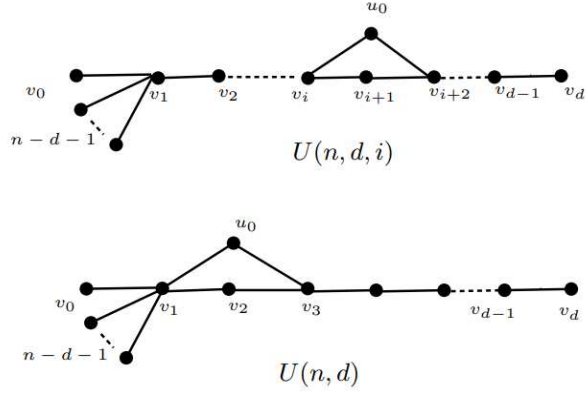}

   \caption{Graph $U(n,d)$ for $4\leq d\leq n-2$.}
    \label{fig:3}
     \end{figure}

Let $U(n,d,i)$(see Fig~\ref{fig:3}) be a graph which consists of the path $v_{1}\cdots v_{d}$ and the cycle $v_{i}u_{0}v_{i+2}v_{i+1}v_{i}$, then attach $n-d-1$ pendant vertices to $v_{1}$. Denote $U(n,d,1)$ and $U(n,d,d-1)$ by $U(n,d)$. Note $d_{G^*}(v_{1})\geq 2$ and $d_{G^*}(v_{i+3})\geq 1$, we have
\begin{align*}
SO_{\alpha}(U(n,d))=(n-d-1)\left((n-d+1)^2+1 \right)^{\alpha}+2\left( (n-d+1)^2+4\right)^{\alpha}+\gamma,
\end{align*}
where $\gamma = 2\cdot 13^{\alpha}+10^{\alpha}$ if $d=4$ and $\gamma = 3\cdot 13^{\alpha}+ (d-5)8^{\alpha}+5^{\alpha}$ if $d\geq 5$.
\begin{lem}\label{lem8}
$SO_{\alpha}(U(n,d,1))\geq SO_{\alpha}(U(n,d,i)) $, where $2 \leq i \leq d-2$.
\end{lem}
\pf
If $i=2$ or $i=d-2$, by Lemma~\ref{lem5} and Lemma ~\ref{lem6} we have
\begin{align*}
SO_{\alpha}(U(n,d))-SO_{\alpha}(U(n,d,2))=& (d_{G^*}(v_{1})-1)\left(\left(d_{G^*}(v_{1})+1\right)^2+1\right)^{\alpha}+2\left(\left(d_{G^*}(v_{1})+1\right)^2+4\right)^{\alpha}\\
&-(d_{G^*}(v_{1})-1)\left(d_{G^*}^{2}(v_{1})+1\right)^{\alpha}- \left(d_{G^*}^2(v_{1})+9\right)^{\alpha}-13^{\alpha}\\
&+\left(4+d_{G^*}(v_{5})\right)^{\alpha}-\left(9+d_{G^*}(v_{5})\right)^{\alpha}\geq 0.
\end{align*}

If $ 3\leq i \leq d-3$, by Lemma~\ref{lem5} and Lemma ~\ref{lem6} we have
\begin{align*}
SO_{\alpha}(U(n,d))-SO_{\alpha}(U(n,d,i))=& (d_{G^*}(v_{1})-1)\left(\left(d_{G^*}(v_{1})+1\right)^2+1\right)^{\alpha}+2\left(\left(d_{G^*}(v_{1})+1\right)^2+4\right)^{\alpha}\\
&-(d_{G^*}(v_{1})-1)\left(d^2_{G^*}(v_{1})+1\right)^{\alpha}- \left(d_{G^*}^2(v_{1})+4\right)^{\alpha}-13^{\alpha}\\
&+\left(4+d_{G^*}(v_{i+3})\right)^{\alpha}-\left(9+d_{G^*}(v_{i+3})\right)^{\alpha}\geq 13^{\alpha} - 8^{\alpha}>0.
\end{align*}
So we get the conclusion.

\begin{figure}[H]
 \centering
  \includegraphics[width=0.5\textwidth]{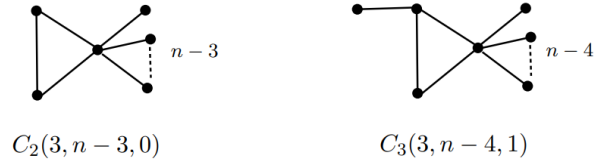}

   \caption{The graphs $C_{2}(3,n-3,0)$ and $C_{3}(3,n-4,1)$}
    \label{fig:2}
     \end{figure}

For $n\geq 4$,  Let $C_{d}(p,q,n-p-q)\in \mathcal{U}_{n,d}$ be a graph obtained from cycle $C_{p}$ by attaching $q$, $n-p-q$ pendent vertices to two neighbor vertices of $C_{p}$, respectively. $\mathcal{U}^{max}_{n,d}$ is the set of unicyclic graphs with maximal general sombor index in $\mathcal{U}_{n,d}$.
\begin{thm}\label{thm3}
If $d=2$, then for $G\in \mathcal{U}_{n,2}$, $SO_{\alpha}(G)\leq SO_{\alpha}(C_{2}(3,n-3,0))$ with equality holds if and only if $G\cong C_{2}(3,n-3,0)$(see Fig~\ref{fig:2}).
If $d=3$, then for $G\in \mathcal{U}_{n,3}$, $SO_{\alpha}(G)\leq SO_{\alpha}(C_{3}(3,n-4,1))$ with equality holds if and only if $G\cong C_{3}(3,n-4,1)$(see Fig~\ref{fig:2}).
\end{thm}

\pf If $G\in \mathcal{U}^{max}_{n,2}$, $g(G)\leq 5$. When $g(G)=3$, then $ G\cong C_{2}(3,n-3,0)$, clearly.  When $g(G)=4$, $G\cong C_{2}(4,n-4,0)$ or $G\cong C_{2}(3,n-4,1)$. If $G\cong C_{2}(4,n-4,0)$, we denote the cycle of $G$ by $C= w_{1}w_{2}w_{3}w_{4}w_{1}$ and $w_{1}$ is connected with $n-4$ pendant vertices. Use the transformation $G\stackrel{w_{1}w{2}}{\longrightarrow}G'$, then $SO_{\alpha}(G)<SO_{\alpha}(G')$ which is a contradiction. If $g(G)=5$, $G$ is $C_{5}$, by Lemma~\ref{lem3} we can get a contradiction. So we get if $G\in \mathcal{U}^{max}_{n,2}$, $G\cong C_{2}(3,n-3,0)$.

If $G\in \mathcal{U}^{max}_{n,3}$, $g(G)\leq 6$. When $g(G)=3$, $G\cong C_{3}(3,n-4,1)$. When $g(G)=4$, $G\cong C_{3}(4,n-4,0)$. Denote the cycle of $G$ by $w_{1}w_{2}w_{3}w_{4}w_{1}$ and $w_{1}$ is connected with $n-4$  pendant vertices. Use the transformation $G\stackrel{w_{3}w_{4}}{\longrightarrow}$, then $G'\in \mathcal{U}_{n,3}$ and $SO_{\alpha}(G)<SO_{\alpha}(G')$, a contradiction. When $g(G)=5$, $G\cong C_{3}(5,n-5,0)$.  Denote the cycle of $G$ by $w_{1}w_{2}w_{3}w_{4}w_{5}w_{1}$ and $w_{1}$ is connected with $n-5$  pendant vertices. Use the transformation $G\stackrel{w_{3}w_{4}}{\longrightarrow}$, then $G'\in \mathcal{U}_{n,3}$ and $SO_{\alpha}(G)<SO_{\alpha}(G')$, a contradiction. If $g(G)=6$, $G$ is $C_{6}$, which cannot have the  maximal $SO_{\alpha}$ index by Lemma~\ref{lem3}. So we get if $G\in \mathcal{U}^{max}_{n,3}$, $G\cong C_{3}(3,n-4,1)$.

\begin{prop}\label{pro2}
Let $G^*\in \mathcal{U}^{max}_{n,d}(3\leq d \leq n-3)$. Then there exist a vertex $v\in VP(G^{*})$ and $G^*-v \in \mathcal{U}_{n-1,d}$.
\end{prop}
\pf Assume to the contrary that for every $v\in VP(G^{*})$, $G^{*}-v \in \mathcal{U}_{n-1,d-1}$. Let $P=diam(G^*)=v_{0}v_{1} \cdots v_{d}$.  By Lemma~\ref{lem3}, we know that $|VP(G^*)\geq 1|$ , suppose $v_{0}\in VP(G^*)$. Note that only $v_{0}$ or $v_{d} \in VP(G^*)$, if there exist other vertex $v' \in VP(G^*)$,  then $G^* - v' \in \mathcal{U}_{n-1,d}$.  By Lemma~\ref{lem4}, for the cycle $C$ in $G^*$, we have $|V(P)\cap V(C)|>1$.

Denote $C= v_{i}v_{i+1}\cdots v_{i+j}u_{l-1}u_{l-2} \cdots u_{0}$, where $u_{0}=v_{i}$. According to Lemma~\ref{lem4}, $j\geq 1$ and $u_{l-1}, u_{l-1}, \cdots, u_{0}$ are not in $P$.  Obviously, $j\leq l$. Since $n\geq d+3$, then $l\geq 3$. We discuss the following two cases.

\textbf{Case 1.} If $j<l$,  we use the transformation $G^{*}\stackrel{v_{i}u_{1}}{\longrightarrow}G^{*'}$,  then $SO_{\alpha}(G^{*})<SO_{\alpha}(G^{*'})$ by Lemma~\ref{lem2},  which is a contradiction.

\textbf{Case 2.} If $j=l$, let $G^{*'}= G^* + v_{i+1}u_{1}+ v_{i+1}u_{2}-v_{i}u_{1}-u_{1}u_{2}$, then by Lemma~\ref{lem5}
\begin{align*}
SO_{\alpha}(G^*)-SO_{\alpha}(G^{*'})=&
\left(d_{G^*}^2(v_{i-1})+9\right)^\alpha-\left(d_{G^*}^2(v_{i-1})+4\right)^\alpha+2\cdot13^\alpha+ 2\cdot8^\alpha- 3\cdot 20^\alpha -17^\alpha \\
<& 10^\alpha-5^\alpha + 2 \cdot 13^\alpha+ 2\cdot8^\alpha-3\cdot 20^\alpha-17^\alpha<0.
\end{align*}
In both cases, $SO_{\alpha}(G^*)<SO_{\alpha}(G^{*'})$. So $G^*\notin \mathcal{U}^{max}_{n,d}$, we get a contradiction. \qed
\begin{prop}\label{pro3}
Let $G^* \in \mathcal{U}_{n,d}$. Denote $VP^{*}(G^*) = \left\{v|v\in VP(G^*), G^*-v \in \mathcal{U}_{n-1,d}\right\}$.
Let $u\in \bigcup_{v\in VP^{*}}$, $NP(u)=\left\{w\in N_{G^*}(u)| \text{w is not a pendant vertex}\right\}$.
Then there exist a vertex $v_{0}\in VP(G^{*})$ and $G^*-v_{0} \in \mathcal{U}_{n-1,d}$, for $u_{0}\in N_{G^*}(v_{0})$, $NP(u_{0})\geq 2$.
\end{prop}

According to Proposition~\ref{pro2}, $VP^*(G^*)\geq 1$. Assume to the contrary that for all $u \in N_{G^*}(v)$(where $v\in VP^{*}(G^*)$), $NP(u)=1$.

\textbf{Claim 1.} $VP^*(G^*)\subseteq N_{G^*}(v_{1})\cup N_{G^*}(v_{d-1})$. Suppose that there exist a vertex $u_{0} \in VP^*(G)$ and $u_{0} \notin N_{G^*}(v_{1})\cup N_{G^*}(v_{d-1})$. $N_{G^*}(u_{0}) = u_{1}$, according to the assumption, we have $NP(u_{1}) = 1$. Note that $u_{1}\notin V(C)\cup V(P)$, otherwise $NP(u_{1}) = 2$. Assume $u_{1}u_{2}$ is a non-pendant edge, we use the transformation $G^{*}\stackrel{u_{2}u_{1}}{\longrightarrow}G^{*'}$, then $SO_{\alpha}(G^*)<SO_{\alpha}(G^{*'})$ according to Lemma~\ref{lem2}. So Claim 1 holds.

 We denote $C= v_{i}v_{i+1}\cdots v_{i+j}u_{l-1}u_{l-2} \cdots u_{0}$, where $u_{0}=v_{i}$ and $j\geq 1$ according to Lemma~\ref{lem4}. Note that $\left\{u_{1}, u_{2},\cdots u_{l-1} \right\} \nsubseteq V(P)$, $l\geq 2$. It's clearly that $j\leq l$.

\textbf{Claim 2.} $|V(C)\setminus V(P)|=1$.  Assume to the contrary that $|V(C)\setminus V(P)|\geq 2$,  we discuss the following two cases.

\textbf{Case 1.} If $j<l$,  apply the transformation $G^{*}\stackrel{v_{i}u_{1}}{\longrightarrow}G^{*'}$, then $SO_{\alpha}(G^{*})<SO_{\alpha}(G^{*'})$ by Lemma~\ref{lem2}, which is a contradiction.

\textbf{Case 2.} If $j=l$, let $G^{*'}=G^{*}-v_{i}u_{1}-u_{1}u_{2}+v_{i+1}u_{1}+v_{i+1}u_{2}$, then $G^{*'}\in \mathcal{U}_{n,d}$. We have $d_{G^*}(v_{i-1})\geq 1$, then by Lemma~\ref{lem5}
\begin{align*}
SO_{\alpha}(G^*)-SO_{\alpha}(G^{*'})&= \left(d_{G^*}^2(v_{i-1})+9\right)^{\alpha}-\left(d_{G*}^2(v_{i-1})+4\right)^{\alpha}+2\cdot13^{\alpha}+2\cdot 8^{\alpha}-3\cdot20^{\alpha}-17^{\alpha}\\
&\leq 10^{\alpha}-5^{\alpha}+2\cdot13^{\alpha}+2\cdot 8^{\alpha}-3\cdot20^{\alpha}-17^{\alpha}<0,
\end{align*}
when $0<\alpha<1$. So we get a contradiction.

\textbf{Claim 3.} $d_{G^*}(v_{d})=1$. On the contrary, we suppose $d_{G^*}(v_{d})=2$. Note that $N(v_{d-1})\cap VP(G^*) = {\O}$, otherwise $|NP(v_{d-1})|\geq 2$. Note that $d_{G^*}(v_{d-3})\geq 2$. By Claim 2, $|V(C)\setminus V(P)| = 1$, thus we consider two cases.

\textbf{Case 1.} $C= v_{d-1}v_{d}u_{1}v_{d-1}$. Let $G^{*'}= G^{*}+u_{1}v_{d-1}+u_{1}v_{d-2}-u_{1}v_{d} $, we have $G^{*'}\in \mathcal{U}_{n,d}$ and
\begin{align*}
\left(d_{G^*}^2(v_{d-3})+4\right)^{\alpha} -\left( d_{G^*}^2(v_{d-3})+9\right)^{\alpha} + 13^{\alpha}-18^{\alpha}+8^{\alpha}-10^{\alpha}<0.
\end{align*}

\textbf{Case 2.} $C= v_{d-2}v_{d-1}v_{d}u_{1}v_{d-2}$. Let $G^{*'} = G^{*} + u_{1}v_{d-1}- u_{1}v_{d}$, then $G^{*'}\in \mathcal{U}_{n,d}$ and
\begin{align*}
2(8^{\alpha}-10^{\alpha})<0.
\end{align*}
By the above two cases, we get $G^{*}$  is not the graph with maximal $SO_{\alpha}$ index in $\mathcal{U}_{n,d}$.

\textbf{Claim 4.} $|E(C)| = 4$. Assume to the contrary that $|E(C)|=3$. Denote $C= v_{i}u_{1}v_{i+1}v_{i}(2\leq i\leq d-2)$. Note that $d_{G^*}(v_{1})\geq 3$ and $d_{G^*}(v_{i+2})\geq 1$.

\textbf{Case 1.} $i=2$. Let $G^{*'}= G^{*}+v_{1}u_{1}-v_{2}u_{1}$, then $G^{*'}\in \mathcal{U}_{n,d}$ and
\begin{align*}
SO_{\alpha}(G^*)-SO_{\alpha}(G^{*'})=&(d_{G^*}(v_{1})-1)\left( \left(1+d_{G^*}^2(v_{1})\right)^{\alpha}-\left(1+(d_{G^*}(v_{1}+1))\right)^{\alpha}\right)+18^{\alpha}\\
&-\left((d_{G^*}(v_{1}+1))^2+4\right)^{\alpha}
\leq 18^{\alpha}-20^{\alpha}<0.
\end{align*}

\textbf{Case 2.} $i\geq 3$. Let $G^{*'} = G^{*}-u_{1}v_{i}- u_{1}v_{i+1}+u_{1}v_{1}+u_{1}v_{3}$, then $G^{*'}\in \mathcal{U}_{n,d}$ and by Lemma~\ref{lem5} and Lemma~\ref{lem6}, we have
\begin{align*}
SO_{\alpha}(G^{*})-SO_{\alpha}(G^{*'})=&(d_{G^*}(v_{1})-1)\left(\left(d_{G^*}^2(v_{1})+1 \right)^{\alpha}-\left( \left(d_{G^*}(v_{1})+1\right)^2+1\right)^{\alpha}\right) \\
&+\left( d_{G^*}^2(v_{1})+9\right)^{\alpha}- 2\left( \left(d_{G^*}(v_{1})+1\right)^2+4\right)^{\alpha}+13^{\alpha}\\
&+ \left(d^2_{G^*}(v_{i+2})+9\right)^{\alpha}-\left(d^2_{G^*}(v_{i+2})+4\right)^{\alpha}\\
\leq& 2\left(10^{\alpha}-17^{\alpha}\right)+2\left(13^{\alpha}-20^{\alpha}\right)+ 10^{\alpha}-5^{\alpha}<0.
\end{align*}
So $SO_{\alpha}(G^*)<SO_{\alpha}(G^{*'})$ and $G^* \notin \mathcal{U}^{max}_{n,d}$. Claim 4 is proved.

By the above discussion, we get $G^*\cong U(n,d,i)$, where $ 2\leq i\leq d-3$, since $d_{G^*}(v_{d})=1$(Claim3). By Lemma~\ref{lem8}, we have $SO_{\alpha}(G^*)< SO_{\alpha}(U(n,d)) $, which is a contradiction and the proof is complete.\qed

\begin{thm}\label{thm1}
Let $G\in \mathcal{U}_{d+2,d}(d\geq4)$. Then $SO_{\alpha}(G)\leq SO_{\alpha}(U(d+2,d))$, with equality holds if and only if $G\cong U(d+2,d)$.
\end{thm}

\pf Suppose that $G^* \in \mathcal{U}^{max}_{d+2,d}$. Let $P=diam(G^*)=v_{0}v_{1}\cdots v_{d}$. Since $n=d+2$, there exists a vertex $u\notin V(P)$. Denoted the cycle of $G^*$ by $C=v_{i}v_{i+1}uv_{i}$ or $C=v_{i}v_{i+1}v_{i+2}uv_{i}$.

\textbf{Claim 1.} $\left\{v_{0},v_{d}\right\}\cap V(C)= \O$. Suppose $v_{0}\in V(C)$, then $C=v_{0}v_{1}uv_{0}$ or $C=v_{0}v_{1}v_{2}uv_{0}$.

\textbf{Case 1.} $C=v_{0}v_{1}uv_{0}$. Let $G^{*'}=G'-v_{0}u-v_{1}u+v_{1}u+v_{2}u$, then  $G^{*'}\in \mathcal{U}_{d+2,d}$. We have
$$SO_{\alpha}(G^*)- SO_{\alpha}(G^{*'})= 2\cdot8^{\alpha}-10^{\alpha}-18^{\alpha}<0.$$

\textbf{Case 2.} $C=v_{0}v_{1}v_{2}uv_{0}$. Let $G^{*'}= G^{*}-v_{0}u+v_{1}u$., then $G^{*'} \in \mathcal{U}_{d+2,d}$. We have
$$SO_{\alpha}(G^*)-SO_{\alpha}(G^{*'})=2\cdot8^{\alpha}-10^{\alpha}-18^{\alpha}<0.$$
Both of two cases contradict with the assumption, thus $v_{0}\notin V(C)$.  We can get $v_{d}\notin V(C)$ similarly.

\textbf{Claim 2.} $|E(C)|=4$. If $|E(C)|=3$, denote $C = v_{i}v_{i+1}u_{1}v_{i}(1\leq i\leq d-2)$.  Since $d\geq 4$,  we have $d_{G^*}(v_{i-1})=2$ or $d_{G^*}(v_{i+2})=2$.  Suppose $d_{G^*}(v_{i+2})=2$, then $i\leq d-3$ and $d_{G^*}(v_{i+3})=1$ or $2$. Let $G^{*'}= G^{*} + u_{1}v_{i+2}-u_{1}v_{i+1}$, then we have $G^{*'}\in \mathcal{U}_{n,d}$ and
\begin{align*}
SO_{\alpha}(G^*)-SO_{\alpha}(G^{*'})=18^{\alpha}-13^{\alpha}+\left(d_{G^*}^2(v_{i+3})+4\right)^{\alpha}-\left(d_{G^*}^2(v_{i+3})+9\right)^{\alpha}<0,
\end{align*}
since if $d_{G^*}(v_{i+3})=1$, $18^{\alpha}-13^{\alpha}+5^{\alpha}-10^{\alpha}<0$ when $0<\alpha<1$. And if $d_{G^*}(v_{i+3})=2$,  $18^{\alpha}+8^{\alpha}-2\cdot 13^{\alpha}<0$ when $0<\alpha<1$. So we get $SO(G^{*})<SO(G^{*'})$, which is a contradiction.

So by Claim 1 and Claim 2, we know $G^*\cong U(d+2,d,i)$. By Lemma~\ref{lem8}, we get $SO_{\alpha}(G^*)\leq SO_{\alpha}(U(d+2,d))$, which means $G^*\notin \mathcal{U}^{max}_{d+2,d}$. We get a contradiction. The proof is complete.\qed
\begin{thm}\label{thm2}
For any $G\in \mathcal{U}_{n,d}(4\leq d \leq n-2)$, $SO_{\alpha}(G)\leq SO_{\alpha}(U(n,d))$, with equality holds iff $G\cong U(n,d)$.
\end{thm}

\pf If $n= d+2$, conclusion holds according to Theorem~\ref{thm1}. Assume the conclusion holds for $G\in \mathcal{U}_{n-1,d}$.

Let $G^*\in \mathcal{U}^{max}_{n,d}(4\leq d \leq  n-3)$. By Proposition~\ref{pro2}, exist a vertex $v^{*}\in VP(G^{*})$, $G^*-v^*\in \mathcal{U}_{n,d}$, and the neighbor $u$ of $v^{*}$ is adjacent to at least two non-pendant vertices. Denote $N_{G^*}(u)=\left\{v,u_{1},u_{2}\cdots u_{d_{G^*}(u)-1} \right\}$,  we have $d_{G^*}(u_{1})\geq 2$, $d_{G^*}(u_{2})\geq 2$ and $ 3\leq d_{G^*}(u)\leq n-d+1$.

Let $G^{*'} = G^{*} - v^*$, then by Lemma~\ref{lem7}, we have
\begin{align*}
SO_{\alpha}(G^*) =& SO_{\alpha}(G^{*'}) + \left( d^{2}_{G^{*}}(u)+1\right)^{\alpha}
+ \sum_{i=1}^{d_{G^*}(u)} \left( \left(d^2_{G^*}(u)+ d^2_{G^*}(u_{i})\right)^{\alpha}- \left( \left(d^2_{G^*}(u)-1\right)^2+ d^2_{G^*}(u_{i})\right)^{\alpha}\right)\\
\leq & SO_{\alpha}(U(n,d-1)) + \left( d^{2}_{G^{*}}(u)+1\right)^{\alpha}
+ 2\left( \left( d^{2}_{G^*}(u)+4\right)^{\alpha} - \left( \left(d_{G^*}(u)-1\right)^2+4 \right)^{\alpha}\right)\\
&+ \left(d_{G^*}(u)-3\right)\left( \left( d^{2}_{G^*}(u)+1\right)^{\alpha} - \left( \left(d_{G^*}(u)-1\right)^2+1 \right)^{\alpha} \right)\\
=& (n-d-2)\left( \left(n-d \right)^2 + 1\right)^{\alpha} + 2\left(\left(n-d \right)^2+4\right)^{\alpha}+ \gamma + \left(d_{G^*}-2\right)\left( d^2_{G^*}+1\right)^{\alpha} + 2\left(d^2_{G^*}+4\right)^{\alpha} \\
&- 2\left(\left(d_{G^*}-1\right)^2+4\right)^{\alpha}-\left(d_{G^*}(u)-3\right)\left( \left(d_{G^*}(u)-1\right)^2+1 \right)^{\alpha}\\
\leq&  (n-d-2)\left( \left(n-d \right)^2 + 1\right)^{\alpha} + 2\left(\left(n-d \right)^2+4\right)^{\alpha}+ \gamma + (n-d-1)\left( \left(n-d+1\right)^2+1\right)^{\alpha} \\
&+ 2\left((n-d+1)^2+4\right)^{\alpha}-2\left(\left(n-d\right)^2+4\right)^{\alpha}
-(n-d-2)\left((n-d)^2+1\right)^{\alpha}\\
=&SO_{\alpha}(U(n,d)).
\end{align*}
Thus $SO_{\alpha}(G)\leq SO_{\alpha}(U(n,d))$, with equality if and only if $G\cong U(n,d)$.\qed

\end{document}